\documentclass[doublespacing]{elsart}
\begin{document}

\begin{frontmatter}

\title{The property of the set of the real numbers generated by a Gelfond-Schneider operator and the countability of all real numbers}
\author[Slavica]{Slavica Vlahovic}
and
\author[Branislav]{Branislav Vlahovic\corauthref{cor}}
\ead{vlahovic@nccu.edu}
\corauth[cor]{Corresponding author.}

\address[Slavica]{Gunduliceva 2, 44000 Sisak, Croatia}
\address[Branislav]{North Carolina Central University, Durham, NC 27707, USA}

\begin{abstract}



Considered will be properties of the set of real numbers $\Re$ generated by an operator that has form of an exponential function of Gelfond-Schneider type with rational arguments.
It will be shown that such created set has cardinal number equal to ${\aleph_0}^{\aleph_0}=c$. It will be also shown that the same set is countable. The implication of this contradiction to the countability of the set of real numbers will be discussed.

\end{abstract}


\end{frontmatter}

\section{Introduction}

In 1900 D. Hilbert announced a list of twenty-three outstanding unsolved problems. The seventh problem was settled in 1934 by A. O. Gelfond and an independent proof by Th. Schneider in 1935.  They proved that if $\alpha$ and $\beta$ are algebraic numbers with $\alpha \neq 0, \alpha \neq 1$, and if $\beta$ is not a real rational number, then any value of ${\alpha}^{\beta}$ is transcendental [1, 2].

For instance transcendental number is $2^{\sqrt{2}}$.

This can be written in the form of
\begin{equation}
{({2\over1})}^{[({2\over1})^{1\over2}]}={({m_{i1}\over n_{i1}})}^{[({{m_{i2}\over n_{i2}})}^{({m_{i3}\over n_{i3}})}]}
\label{gelement}
\end{equation}
where $m_i, n_i \in N$.

We can ask ourselves a following question which "class" of transcendental numbers can be presented this way? Or can any transcendental number be expressed in the form (\ref{gelement}). Answer is obvious; some transcendental numbers cannot be expressed this way, for instance number $e$ cannot be presented by
\begin{equation}
e={m_1\over n_1}^{[{m_2\over n_2}^{m_3\over n_3}]}
 \label{e}
    \end{equation}
because after taking logarithm from both sides
\begin{equation}
1={m_2\over n_2}^{m_3 \over n_3}ln{m_1\over n_1},
 \label{le}
    \end{equation}
and this cannot be, because $ln{m_1\over n_1}$ is always transcendental [3-5] for $m_1,n_1\in N$.

However, one can take more freedom and try to express the number $e$ in the form
\begin{equation}
{[{m_1\over n_1}^{m_2\over n_2}]}^{[{m_3\over n_3}^{m_4\over n_4}]}
 \label{ec}
    \end{equation}

or even more freedom and try to present the number $e$ in the form
${a_1}^{a_2}$, where both ${a_1}$ and ${a_2}$ can have the form (\ref{ec}).
Obviously, the argument such as shown in (\ref{le}), that number $e$ cannot be presented in such way, cannot be applied anymore since both ${a_1}$ and ${a_2}$ can be now transcendental numbers.

One can go even further (as it is done in [6]) and take much more freedom in generating the numbers or a set of numbers, which elements will be generated through a general element of the sequence that has the form:
\begin{equation}
{a_1}^{{a_2}^{{a_3}^{.^{.^{.^{{a_n}^{.^{.^{.}}}}}}}}}
\label{gsequence}
\end{equation}
where in (\ref{gsequence}) each element $a_i$ of bases and exponents has the
following form:
\begin{equation}
a_i={[({{m_{i1}\over n_{i1}})}^{({m_{i2}\over n_{i2}})}]}^{[({{m_{i3}\over n_{i3}})}^{({m_{i4}\over n_{i4}})}]}
\label{gelementn}
\end{equation}
where $m_{ij},n_{ij}\in N, i=1, 2, 3,...n, j = 1, 2, 3, 4$.

The question remains: which class of the transcendental numbers can be or can not be represented in this way? Can majority of the transcendental numbers be presented or can not be presented in this way? If some transcendental numbers can not be presented, is that set countable or not?

First let us note that the set of numbers generated through the operator (\ref{gsequence}) looks similar to the set of the real numbers.  Such set does not have the first and last element, it has subset of all rational numbers, and it is dense everywhere in rational, algebraic and transcendental numbers. However, it may not be equal to the set of the real numbers since it is harder to prove that it is dense in Dedekind's sense, since this would require proof that it does not have holes, i.e. that all numbers can be represented in this way.

To avoid that difficulty, let as assume that some numbers can not be presented in this way and let us focus here only on estimating the number of the elements in such set, i.e. on determining the cardinality of such set of numbers.

\section{The cardinality of the generated set of numbers}

Let us generate the set of the real numbers through relation (\ref{gsequence}) where each base and exponent element $a_i$ has the form (\ref{gelementn}). The mechanism to generate the elements of the set is to write (\ref{gsequence}) for all possible combinations of arguments, with the sum of all bases and exponents equal to 2, 3, 4,... and so on.  As the sum increases the number of the exponents will expand. The sample of such generated set with a procedure to avoid double counting of the same numbers is given in [6]. However, let us focus here on our main task which is to estimate the cardinality of such generated set when the process described above continues to infinity.

For each particular number the general element (\ref{gsequence}) that corresponds to that number will have a final number of the exponents $a_i$. However, since the process of generating new numbers continues to infinity there is no an upper limit for the number of the exponents $a_i$ that will be generated by the general element (\ref{gsequence}), which will also go to infinity. Each of the elements of $a_i$ will have $\aleph_0$ possible combinations. This is obvious, since for any arbitrary large value $n \in N$ which one could take for the number of combinations, that value will be exceeded in this described process. The same is true for the number of the exponents. The number of the exponents will also be $\aleph_0$, since again any arbitrary taken number that one could chose for the value of the number of the exponents (does not matter how large is the number) will be exceeded in the described process, which continues to infinity.

Therefore, the above described set will have ${\aleph_0}^{\aleph_0} = c$ elements which makes it equivalent in the cardinal number to the set of the real numbers. A one to one correspondence between such produced set and the set of natural numbers $N$ can be easily obtained by arranging the set elements by the sum of the exponents, as it is done for instance in [6].

We will not here proceed to discuss what could be wrong with the Cantor's famous diagonal proof of countability of the set of real numbers; some of the relevant remarks are done in [6, 7]. Let us note here that that proof could be wrong since it uses the method of induction which, as it is well known [8, 9, 10], can not be applied on the infinite sets. With that method one can only prove that a number created by the diagonal procedure can be different from any $n$ numbers in the set. The method can not prove that that number is different from any number in the assumed denumerable set, which has infinite number of the elements. So, one can move through that set using the diagonal procedure to higher positions numbers $n$ in the sequence, but can not go through all the set elements. At least it cannot be done by using the induction method.

\section{Conclusion}

It is proven that the set generated by the general element (\ref{gsequence}) has cardinal number equal to ${\aleph_0}^{\aleph_0} = c$. The same set is also denumerable, the elements can be ordered by the sum of the bases and exponents in (\ref{gsequence}). Therefore it is proven that the cardinality of the real and natural set of numbers are the same, i.e. that ${\aleph_0}^{\aleph_0}= c= \aleph_0$.


\begin{thebibliography}{10}

\bibitem {1} A. O. Gelfond, Doklady Akad. Nauk. S.S.S.R., {\bf 2} (1934), 1-6.
\bibitem {2}Th. Schneider, J. Reine angew. Math., {\bf 172} (1035) 65-69.
\bibitem{3} A. Baker, Linear forms in the logarithms of algebraic numbers I,
II, III, IV, Mathematika, {\bf 13}(1966),204-216; {\bf 14}(1967),102-107, 220-228;
{\bf 15}(1968),204-216.
\bibitem{4} Ch. Hermite, Sur la fonction exponentialle, Oeuvres III, 150-181.
\bibitem{5} A. Baker and D. W. Masser, Transcendence Theory: Advances and
Applications, Academic Press London New York San Francisco, 1977.
\bibitem{6} S. Vlahovic and B. Vlahovic, Countability of the Real Numbers arXiv:math/0403169v1
\bibitem{7} S. Vlahovic and B. Vlahovic, Remarks on Cantor's diagonalization proof of 1891, arXiv:math/0403288
\bibitem{8} A.A. Fraenkel and Y. Bar-Hillel, Foundation of Set Theory,
Amsterdam 1958, chapter IV.
\bibitem{9} A. Frankel, Y. Bar-Hillel and A. Levy, Foundation of set theory,
North Holland, Amsterdam 1973, van Dalen's remarks p. 268.
\bibitem{10} M.
Hallett, Cantorian Set Theory and Limitation of Size, Clarendon
Press Oxford, 1984.


\end{thebibliography}
  	  	\end{document}